\documentstyle[amssymb]{amsart}

\newcommand{\Proof}{\noindent {\sc Proof} \hspace{0.2in}} 

\newtheorem{theorem}{Theorem}
\newtheorem{claim}{Claim}[section]
\newtheorem{lemma}[claim]{Lemma} 
\newtheorem{proposition}[claim]{Proposition}

\newtheorem{definition}[claim]{Definition}

\newcommand{\cf}{{\rm cf}}

\newcommand{\dom}{\operatorname{dom}}
\newcommand{\cl}{\operatorname{cl}}

\newcommand{\lesdot}{\mathrel{\mathord{<}\!\!\raise 0.8
pt\hbox{$\scriptstyle\circ$}}}  

\newcommand{\conc}{^\frown\!}

\newcount\skewfactor
\def\mathunderaccent#1#2 {\let\theaccent#1\skewfactor#2
\mathpalette\putaccentunder}
\def\putaccentunder#1#2{\oalign{$#1#2$\crcr\hidewidth
\vbox to.2ex{\hbox{$#1\skew\skewfactor\theaccent{}$}\vss}\hidewidth}}

\begin{document}

\title[Topological density]{Topological density of ccc Boolean algebras - 
every cardinality occurs}
\author{Mariusz Rabus}
\address{Department of  of Mathematics \\
York University\\
North York, ON, M3J 1P3, Canada}
\author{Saharon Shelah}
\address{Institute of Mathematics \\
The Hebrew University of Jerusalem\\
Jerusalem 91904, Israel}
\thanks{Research supported by Basic Research Foundation, 
administered by the Israeli Academy of
Sciences and Humanities. Publication No. 631}
\email{rabus@@mathstat.yorku.ca $\quad$
shelah@@math.huji.ac.il}
\maketitle
\begin{abstract}
For every uncountable cardinal $\mu$ there is a ccc Boolean 
algebra whose topological density is $\mu$.
\end{abstract}

\section{Introduction}

For a Boolean algebra $B$, the topological density 
$d(B)$ is the minimal cardinal $\mu$ such that there is 
a family $\{D_{\xi}:\xi \in \mu\}$ of ultrafilters of 
$B$ with the union $B\setminus \{0\}$.
Note that if $St(B)$ is the Stone space of $B$, then $d(B)$ is
the density of $St(B)$ (as a topological space).
A Boolean algebra $B$ has the countable chain condition
(ccc) if there is no uncountable collection of pairwise 
disjoint elements of $B\setminus \{0\}$.

The question we consider in this paper is: what cardinals
are topological densities of ccc Boolean algebras?
Hajnal, Juh\'asz, and Szentmikl\'ossy \cite{HJS} prove 
that under some mild set-theoretic assumptions every 
uncountable cardinal is  the topological density of 
some ccc Boolean algebra.
We prove here that the above statement is a theorem  of ZFC.

\begin{theorem}
For every uncountable cardinal $\mu$ there is a ccc Boolean 
algebra $\cal{B}$, such that $d(\cal{B})=\mu$.
\end{theorem}

The rest of the paper is devoted to the proof of the theorem. Let
$\mu$ be an uncountable cardinal. The idea of the proof is to define
$\cal{B}$ as a quotient of a free Boolean algebra generated by 
$\{x_{\nu}:\nu\in T\}$, where $T$ is a set of cardinality $2^{\mu}$.
The reason we index the generators by a set $T$, rather than $2^{\mu}$,
is that an additional structure on $T$ is helpful in defining the 
quotient. In particular, the quotient is defined by imposing a set of 
restrictions of the form $x_{\nu_{0}}\cap x_{\nu_{1}} \cap 
(-(x_{\eta_{0}} \vartriangle x_{\eta_{1}}))=0$, for some 
$\nu_{0},\nu_{1},\eta_{0},\eta_{1}$ in $T$. The definition of $T$
is quite technical but it is the key element of the proof that 
the topological density of our algebra is $\geq \mu$, Lemma 4.1.
The construction of $T$ is done in section 
2. In section 3 we give the definition of the algebra $\cal{B}$
and prove that it has the ccc. In the last section we prove that 
the topological density of $\cal{B}$ is exactly $\mu$.

\section{Preliminaries}

In this section we define the set $T$ and the set of quadruples used in 
the definition of $\cal{B}$. 
For a cardinal $\sigma<\mu$ let $h_{\sigma}:[\sigma^{+}]^{2}
\rightarrow \sigma^{+}$ be such that 
\begin{enumerate}
\item[(a)] it is one-to-one,
\item[(b)] for every $X\in [\sigma^{+}]^{\sigma^{+}}$ and 
$j<\sigma$ there is $j_{1}\in(j,\sigma)$ and $i_{0},i_{1}
\in X$ such that $h_{\sigma}(i_{0},i_{1})=j_{1}$ mod $\sigma$.
\end{enumerate}

To prove that such a function exists it is enough to show that
there is a function $h'_{\sigma}:[\sigma^{+}]^{2}
\rightarrow \sigma$ such that (b) holds. Indeed, then we can define 
$h_{\sigma}$ to be any $1-1$ function from $\sigma^{+}$ to 
$\sigma^{+}$ such that $h_{\sigma}(\alpha,\beta)=
 h'_{\sigma}(\alpha,\beta) \text{ mod } \sigma$.

In order to define $h'_{\sigma}$, first fix  $1-1$ functions 
$q_{\alpha}:\alpha \rightarrow \sigma$ for every $\alpha < \sigma^{+}$.
Now for $\alpha < \beta$  define $h'_{\sigma}(\alpha,\beta)=
q_{\beta}(\alpha)$. To prove (b) let $X\in [\sigma^{+}]^{\sigma^{+}}$.
Let $\beta \in X$ be such that the set $X_{\beta}=
\{\alpha \in X: \alpha < \beta\}$ has cardinality $\sigma$.
Since $q_{\beta}$ is $1-1$ it follows that the image of $[X_{\beta}]^{2}$ 
under $h'_{\sigma}$ is cofinal in $\sigma$, hence (b) holds.

\begin{definition} We define, by induction on $\alpha \leq \mu^{+}$,
$\cal{T}_{\alpha}=(T_{\alpha},\sigma_{\alpha},
P^{l}_{\alpha}, F^{m}_{\alpha})$, $(l<8,m<2)$ such that: 
\begin{enumerate}
\item[(1)] $T_{\alpha}$ is a set of finite sequences,
\item[(2)] $\langle P^{l}_{\alpha}:l<8 \rangle$ is a 
partition of $T_{\alpha}$,
\item[(3)] $\sigma_{\alpha} $ is a function from $T_{\alpha}$
to $\{\sigma: \aleph_{0} \leq \sigma \leq \mu\}$,
\item[(4)] $F^{m}_{\alpha}$ is a partial, two-place, symmetric
function from $P^{7}_{\alpha}$ to $T_{\alpha}$,
\item[(5)] $\cal{T}_{\alpha}$ is increasing, continuous in 
$\alpha$, i.e., if $\beta < \alpha$, then 
\begin{enumerate}
\item[(a)] $T_{\beta} \subseteq T_{\alpha}$, 
\item[(b)] $P^{l}_{\beta}=
P^{l}_{\alpha} \cap T_{\beta}$ for $l<8$, 
\item[(c)] $\sigma_{\beta}=
\sigma_{\alpha} \restriction T_{\beta}$,
\item[(d)] $F^{m}_{\beta}=
F^{m}_{\alpha}\restriction [T_{\beta}]^{2}$ for $m<2$,
\item[(e)] if $\alpha$ is a limit, then 
$T_{\alpha}=\bigcup_{\beta<\alpha}
T_{\beta}$.
\end{enumerate}
\end{enumerate}
\end{definition}

\noindent {\it Case 1.} $\alpha =0$. Define $T_{0}=\{\langle \rangle \}$,
$P_{0}^{0}=\{\langle \rangle \}$, (so $P_{0}^{l}=\emptyset$ for $l\not = 0$).
$\sigma_{0}(\langle \rangle )=\mu$.

\noindent {\it Case 2.} $\alpha =1$. Define $T_{1}=T_{0}\cup
\{\langle \sigma \rangle :\aleph_{0} \leq \sigma \leq \mu\}$,
$P_{1}^{1}=\{\langle \sigma \rangle :\aleph_{0} \leq \sigma \leq \mu\}$,
$\sigma_{1}(\langle \sigma \rangle)=\sigma$.

\noindent {\it Case 3.} $\alpha$ is a limit. Put $T_{\alpha}=
\bigcup_{\beta<\alpha}T_{\beta}$, $P_{\alpha}^{l}=\bigcup_{\beta<\alpha}
P_{\beta}^{l}$, $\sigma_{\alpha}=\bigcup_{\beta<\alpha}\sigma_{\beta}$.

\noindent {\em Case 4.} $\alpha=\beta + 1$. Define $T_{\alpha}=T_{\beta}
\cup\{T_{\alpha,l}:l<8\}$, where

\begin{align*}
T_{\alpha,0}&=\emptyset, \\
T_{\alpha,1}&=\{\eta\conc \langle \sigma \rangle : \eta \in T_{\beta}
\setminus (P^{5}_{\beta}\cup P^{7}_{\beta}), 
\eta\conc \langle \sigma \rangle \not \in T_{\beta},
\aleph_{0} \leq \sigma \leq \sigma_{\beta}(\eta)\}, \\
T_{\alpha,2}&=\{\eta\conc \langle 0 \rangle:\eta \in T_{\beta}
\setminus (P^{5}_{\beta}\cup P^{7}_{\beta}), 
\eta\conc \langle 0 \rangle \not \in T_{\beta}\}, \\
T_{\alpha,3}&=\{\eta\conc \langle \rho\conc \langle 0 \rangle \rangle:
\eta\conc \langle \rho \rangle \in T_{\beta}\setminus (P^{5}_{\beta}
\cup P^{7}_{\beta}), 
\eta\conc \langle \rho\conc \langle 0 \rangle \rangle \not \in T_{\beta}\}, \\
\begin{split}
T_{\alpha}^{4,5}&=\{\eta\conc \langle \rho \rangle: \eta \in T_{\beta},\;
\rho \text{ a sequence of limit length and } \\
&\quad\quad\quad (\forall \zeta < \lg(\rho))
(\eta \conc \langle \rho \restriction \zeta \rangle \in T_{\beta}
\text{ but } \eta\conc \langle \rho \rangle \not \in T_{\beta})\},
\end{split} \\ 
T_{\alpha,4}&=\{\eta\conc \langle \rho \rangle\in T_{\alpha}^{4,5} :
\lg(\rho)<\sigma_{\beta}(\eta)^{+}\}, \\
T_{\alpha,5}&=T_{\alpha}^{4,5}\setminus T_{\alpha,4}, 
\end{align*}
\begin{align*}
\begin{split}
T_{\alpha,6}&=\{\eta\conc \langle \rho\conc\langle \nu_{0},\nu_{1}
\rangle \rangle : \eta\conc \langle \rho \rangle \in T_{\beta}\setminus
(P_{\beta}^{5}\cup P^{7}_{\beta}),\\
&\quad\quad\quad \text{ and } \eta\conc \langle \rho \rangle
\vartriangleleft \nu_{l}, l=0,1, \text{ and }
\nu_{0},\nu_{1}\in T_{\beta}, \nu_{0}\not =\nu_{1}\}, 
\end{split} \\
T_{\alpha,7}&=\{\eta\conc \langle \rho \rangle\conc \langle i \rangle :
\eta\conc \langle \rho \rangle \in P_{\beta}^{5}, \mbox{ \rm and }
i<\sigma_{\beta}(\eta)^{+}, \mbox{ \rm and }
\eta\conc \langle \rho \rangle\conc \langle i \rangle \not \in T_{\beta}\}.
\end{align*}

\noindent Let 
$P_{\alpha}^{l}=P_{\beta}^{l} \cup T_{\alpha,l}$ for $l<8$, 
and define $\sigma_{\alpha}$ by: 

\begin{equation*}
\sigma_{\alpha}(\tau)=
\begin{cases}
\sigma_{\beta}(\tau) &\text{if } \tau\in T_{\beta},\\
\tau(n) &\text{if } \tau\in T_{\alpha,1},\lg(\tau)=n+1,\\
\sigma_{\beta}(\tau\restriction n) &\text{if } 
\tau \in \{T_{\alpha,l}:l\not =1\},
\lg(\tau)=n+1.
\end{cases}
\end{equation*}

Finally define $F_{\alpha}^{m}$, $m=0,1$.  $F_{\alpha}^{m}(\tau_{1},\tau_{2})$
is well defined
if $F_{\beta}^{m}(\tau_{1},\tau_{2})$ is well defined, or 
for some  
$\eta\conc \langle \rho\rangle \in P_{\beta}^{5}$  for $l=1,2$, we have
$\tau_{l}=\eta\conc\langle \rho \rangle \conc\langle i_{l}\rangle$,
$i_{1}\not =i_{2}$, and $\rho(\zeta)$ is a pair $(\nu_{0},\nu_{1})$, i.e., 
$\eta\conc \langle \rho\restriction \zeta \conc \langle \nu_{0},\nu_{1}
\rangle\rangle \in P^{6}$,
where $\zeta=h_{\sigma(\eta)}(i_{1},i_{2})$.
In the first case define $F_{\alpha}^{m}(\tau_{1},\tau_{2})=
F_{\beta}^{m}(\tau_{1},\tau_{2})$. In the second case define
$F_{\alpha}^{m}(\tau_{1},\tau_{2})=\nu_{m}$.

\begin{proposition} 
\mbox{\em (1)} $T_{\alpha}$ is well defined for $\alpha\leq \mu^{+}$.
\begin{enumerate}
\item[(2)] Each member of $T_{\alpha}$, ($\alpha\leq \mu^{+}$) 
is a finite sequence.
\item[(3)] For every $\eta \in T_{\mu^{+}}$, the sequence
$\langle \sigma(\eta\restriction k):k\leq \lg(\eta) \rangle$ is non-increasing.
\item[(4)] If $\eta\conc\langle \rho \rangle \in P^{5}$ and 
$m_{1},m_{2} <2$, $\tau_{1},\tau_{2},\tau_{3},\tau_{4} \in
\{\eta\conc\langle \rho \rangle\conc \langle i\rangle : i< \sigma(\eta)^{+}\}$
and $F_{\alpha}^{m_{1}}(\tau_{1},\tau_{2})=
F_{\alpha}^{m_{2}}(\tau_{3},\tau_{4})$,
then $m_{1}=m_{2}$, $\{\tau_{1},\tau_{2}\}=\{\tau_{3},\tau_{4}\}$. Moreover, 
the conclusion holds if we assume that 
$F_{\alpha}^{m_{1}}(\tau_{1},\tau_{2})\restriction \lg(\eta\conc \rho)=
F_{\alpha}^{m_{2}}(\tau_{3},\tau_{4})\restriction \lg(\eta\conc \rho)$. 
\item[(5)] If $\eta\in P_{\mu^{+}}^{7}$,
then $\eta$ is maximal in $(T_{\mu^{+}},\vartriangleleft)$.
\item[(6)] $|T_{\mu^{+}}|=2^{\mu}$.
\end{enumerate}
\end{proposition}

\Proof Straightforward.

\medskip
\noindent Let $T=T_{\mu^{+}}$, $P^{l}=P^{l}_{\mu^{+}}$ for $l<8$, 
$\sigma=\sigma_{\mu^{+}}$ and $F_{m}=F_{\mu^{+}}^{m}$, $m=0,1$.

\begin{definition} \begin{enumerate}
\item[(1)] We say that $X\subseteq T$ is $1$-closed if:
\begin{enumerate}
\item[(a)] $\langle \rangle \in X$,
\item[(b)] if $\eta \vartriangleleft \eta_{1}$, $\eta_{1}\in X$, then 
$\eta \in X$,
\item[(c)] if  $\eta\conc\langle \rho \rangle \in P^{5}$ and for $k=1,2$, 
$\tau_{k}=\eta\conc\langle \rho \rangle\conc\langle i_{k}\rangle \in X$,
$i_{1}\not =i_{2}$ and $m<2$, then $F_{m}(\tau_{1},\tau_{2})
\in X$ if it is well defined.
\end{enumerate}
\item[(2)] We say that $X\subseteq T$ is $2$-closed if it 
is $1$-closed and:
\begin{enumerate}
\item[(d)] if 
$\eta\conc\langle\rho\conc\langle\nu_{0},\nu_{1}\rangle\rangle \in 
P^{6}\cap X$, then $\nu_{0},\nu_{1}\in X$,
\item[(e)] if $\eta\conc\langle \rho_{1} \rangle$,
$\eta\conc\langle \rho_{2} \rangle \in X $, 
$\zeta=\sup\{\xi:\rho_{1}\restriction \xi =\rho_{2}\restriction \xi\}$,
then
$\eta\conc\langle \rho_{1}\restriction (\zeta+1) \rangle \in X$ if 
$\zeta<\lg(\rho_{1})$, and $\eta\conc\langle \rho_{2}\restriction 
(\zeta+1) \rangle \in X$ if 
$\zeta<\lg(\rho_{2})$.
\end{enumerate}
\end{enumerate}
\end{definition}

\begin{proposition}\begin{enumerate}
\item[$(1)$] $T_{\alpha}$ is closed for $\alpha\leq \mu^{+}$.
\item[(2)] The family of $k$-closed sets is closed under
intersections, $k=1,2$. 
\item[(3)] If $X\subseteq T$ is finite, then $\cl_{k}(X)$ is finite,
$k=1,2$.
\end{enumerate}
\end{proposition}

\Proof (1), (2) are straightforward. To prove (3), prove by 
induction on $\alpha$, that if $X\subseteq T_{\alpha}$
is finite, then $\cl_{k}(X)$ is finite.

\section{Definition of the algebra, and the ccc}

In this section we, first, define the algebra, and second,
prove that it has the ccc. The proof is preceded by two 
propositions, which give a sufficient condition for an element 
of the algebra to be non-zero.

\begin{definition} \mbox{\em (1) } $\cal{B}_{T}$ 
is the Boolean algebra generated 
by $\{x_{\eta}:\eta\in T\}$ freely, except the equations in the 
following set:
\begin{multline*}
\Gamma=\{\bold{e}_{\tau_{1},\tau_{2}}=
[x_{\tau_{1}}\cap x_{\tau_{2}} \cap (- (x_{F_{0}(\tau_{1},\tau_{2})} 
\vartriangle
x_{F_{1}(\tau_{1},\tau_{2})}))=0] :\\
\tau_{1},\tau_{2}\in T, \mbox{ \rm and }
F_{m}(\tau_{1},\tau_{2}) \mbox{ \rm is well defined, } m=0,1.\}
\end{multline*}

\noindent \mbox{\em (2)} For $X \subseteq T$ let 
\begin{equation*}
\Gamma_{X}=\{
{\bold{e}}_{\tau_{1},\tau_{2}}:\tau_{1},\tau_{2}\in X, \text{ and }
F_{l}(\tau_{1},
\tau_{2}) \text{ is well defined }l=0,1. \} 
\end{equation*}
 
\noindent \mbox{\em (3)} For $\alpha < \mu$ define $\cal{B}_{T_{\alpha}}$
to be the subalgebra of $\cal {B}_{T}$ generated by $\{x_{\eta}:
\eta\in T_{\alpha}\}$.

\noindent \mbox{\em (4)} $\cal{B}_{0}$ is the trivial Boolean algebra with 
the universe $\{0,1\}$.
\end{definition}

Note: for $\eta\in T$ we consider $x_{\eta}$ to be an element of
${\cal B}_{T}$, i.e., 
it is an equivalence class of the element $x_{\eta}$.

\begin{proposition}\label{p1}
For a Boolean term $\bf{t}=t(y_{0},\ldots , y_{n-1})$ and 
$\eta_{0}, \ldots ,\eta_{n-1} \in T$, $\cal{B}_{T} \vDash t(x_{\eta_{0}},
\ldots ,x_{\eta_{n-1}})>0$ if and only if
there is a function $f:T\rightarrow \{0,1\}$ such that 
$\cal{B}_{0}\vDash t(f(\eta_{0}),
\ldots ,f(\eta_{n-1}))=1$ and $(*)_{f,T}$ holds, where for $X\subseteq T$
we define:

\begin{enumerate}
\item[$(*)_{f,X}$]
 If  
${\bold{e}}_{\tau_{1},\tau_{2}} \in \Gamma_{X}$  and 
$f(\tau_{1})=1=f(\tau_{2})$, then  $f(F_{0}(\tau_{1},\tau_{2})) \not =
 f(F_{1}(\tau_{1},\tau_{2}))$.
\end{enumerate}
\end{proposition}

\Proof 
(1) Assume that $f:T\rightarrow \{0,1\}$ is such that 
$\cal{B}_{0}\vDash t(f(\eta_{0}),
\ldots ,f(\eta_{n-1}))=1$ and $(*)_{f,T}$ holds. Note that 
${\cal B}_{T}\models  t(x_{\eta_{0}},
\ldots ,x_{\eta_{n-1}}) >0$ if and only if there is a homomorphism
$h: {\cal B}_{T}\rightarrow {\cal B}_{0}$ such that ${\cal B}_{0}\models
h( t(x_{\eta_{0}},
\ldots ,x_{\eta_{n-1}}))$$=1$.
The function $f:T\rightarrow \{0,1\}$ defines a homomorphism $\bar{f}$
from a free algebra generated by $\{x_{\eta}:\eta\in T\}$ into
${\cal B}_{0}$. Such homomorphism induces an homomorphism of
${\cal B}$ into ${\cal B}_{0}$ if and only if 
$\bar{f}(x_{\tau_{1}}\cap x_{\tau_{2}} \cap (- (x_{F_{0}(\tau_{1},\tau_{2})} 
\vartriangle
x_{F_{1}(\tau_{1},\tau_{2})})))=0$ for $\tau_{1},\tau_{2}$ such that 
$F_{m}(\tau_{1},\tau_{2})$ are well-defined. Clearly this is 
equivalent to $(*)_{f,T}$.

\smallskip
\noindent (2) Assume $\cal{B}_{T} \vDash t(x_{\eta_{0}},
\ldots ,x_{\eta_{n-1}})>0$. Without loss of generality
$t(x_{\eta_{0}},
\ldots ,x_{\eta_{n-1}})$ $=\bigcap_{l<n}x_{\eta_{l}}^{\epsilon(l)}$,
where $\epsilon:n\rightarrow \{0,1\}$, and $x^{1}=x$, and $x^{0}=-x$.
Moreover, we can assume that $\{\eta_{0}, \ldots ,
\eta_{n_1}\}$ is $1$-closed. Define $f:T
\rightarrow \{0,1\}$ by: $f(\eta_{l})=\epsilon(l)$ for $l<n$, and
$f(\rho)=0$ for $\rho \not \in\{\eta_{0}, \ldots ,
\eta_{n_1}\}$. Clearly ${\cal B}_{0} \models t(f(\eta_{0}),\ldots ,
f(\eta_{n-1}))=1$ and $(*)_{f,T}$ holds.

\begin{proposition}\label{p2} \begin{enumerate} \item[(1)] If 
$X\subseteq T $ is $1$-closed, $f:X\rightarrow \{0,1\}$,  
$(*)_{f,X}$ holds,  $\eta_{0},\ldots , \eta_{n-1} \in X$, and 
$t(y_{0}, \ldots ,y_{n-1})$ is a Boolean term such that 
$\cal{B}_{0} \vDash t(f(\eta_{0}), 
\ldots ,f(\eta_{n-1}))=1$, then $\cal{B}_{T} \vDash t(x_{\eta_{0}},
\ldots ,x_{\eta_{n-1}}) >0$.
\item[(2)] For $\eta \not = \nu$ in $T$, $\cal{B}_{T} \vDash x_{\eta}
\not = x_{\nu}$, moreover  $\cal{B}_{T} 
\vDash (x_{\eta}  \setminus x_{\nu}) >0$.
\item[(3)] If $X \subseteq Y \subseteq T$, $\cl_{1}(X)\subseteq Y$, 
$\eta_{0}, \ldots ,\eta_{n-1} \in Y$, 
$f:Y \rightarrow \{0,1\}$, $X \subseteq f^{-1}(\{1\})$, 
$t(y_{0}, \ldots ,y_{n-1})$ is a Boolean term such that

\noindent $\cal{B}_{0}\vDash t(f(\eta_{0}), \ldots ,f(\eta_{n-1}))=1$,
and 
 ${\bold{e}}_{\tau_{1},\tau_{2}} \in \Gamma_{X}$ implies 
that 

\noindent $|\{F_{0}(\tau_{1},\tau_{2}),F_{1}(\tau_{1},\tau_{2})\}
\cap f^{-1}(\{1\})|=1$, then 
$\cal{B}_{T}\vDash t(x_{\eta_{0}}, \ldots ,x_{\eta_{n-1}})>0$.
\end{enumerate}
\end{proposition}

\Proof (1) Use \ref{p1} for $f \cup 0_{T\setminus Y}$.
For (2) use part (1) with $X=T$, $f(\eta)=1$ and 
$f(\rho)=0$ for $\rho\not =\eta$. (3) is similar to (1).

\begin{lemma} $\cal{B}_{T}$ satisfies the ccc., in fact a strong version of
the ccc: for every collection of $\kappa=\cf(\kappa) >\aleph_{0}$ 
elements
of $\cal{B}_{T}$, there is a subcollection of size $\kappa$ which 
generates a filter.
\end{lemma}

\Proof Let $\cal{B}_{T} \vDash a_{\alpha} >0$ for $\alpha <\kappa$. 
Let $a_{\alpha}=t_{\alpha}(x_{\eta_{\alpha,0}}, \ldots ,x_{\eta_{\alpha,
n_{\alpha}-1}} )$, each $t_{\alpha}$ is a Boolean term.
Without loss of generality we can assume that:
\begin{enumerate}
\item[(1)] $\{\eta_{\alpha,0}, 
\ldots ,\eta_{\alpha,n_{\alpha}-1}\}$ is $2$-closed for each $\alpha$,
\item[(2)] $t_{\alpha}=t$, $n_{\alpha}=n$, 
\item[(3)] $\langle \{\eta_{\alpha ,k}:k<n\}: \alpha<\kappa
 \rangle$ is a $\Delta$-system,
i.e., for some $n(*)\leq n$ we have: if $k<n(*), \alpha<\kappa$ then 
$\eta_{\alpha,k}=\eta_{k}$, 
and $\langle \{\eta_{\alpha,k}:k\in[n(*),n)\}:\alpha<\kappa\rangle$ 
is a sequence of pairwise disjoint sets.
\end{enumerate}

We can assume that $n>n(*)$, as otherwise $a_{\alpha}=a_{0}$ for every
$\alpha<\kappa$, and we are done.
Let $f_{\alpha}:T\rightarrow \{0,1\}$ be such that $(*)_{f_{\alpha},T}$
holds and 

\noindent $\cal{B}_{0}\vDash t(f_{\alpha}(\eta_{\alpha,0}), \ldots ,
f_{\alpha}(\eta_{\alpha,n-1}))=1$. Without loss of generality we can assume
that: 
\begin{enumerate}
\item[(4)]
$f_{\alpha}(\eta_{\alpha,k})=\bf{t}_{k}$, i.e.,
does not depend on $\alpha$, 
\item[(5)] the truth values of 
``$\eta_{\alpha,k_{1}} \vartriangleleft \eta_{\alpha,k_{2}}$'',
``$\lg(\eta_{\alpha,k})=m$'',
``$\eta_{\alpha,k}\in P^{l}$'' do not depend on $\alpha$.
\item[(6)]
if $\eta\conc \langle \rho \rangle \in P^{5} \cap \{\eta_{k}:k<n(*)\}$,
$\tau_{1}=\eta\conc \langle \rho \rangle\conc\langle i_{1} \rangle =
\eta_{\alpha,k_{1}},$ $k_{1}\in [n(*),n)$, and 
$\tau_{2}=\eta\conc \langle \rho \rangle\conc\langle i_{2} \rangle =
\eta_{\alpha,k_{2}},$ $k_{2}\in [n(*),n)$, $i_{1}\not =i_{2}$, 
{\em then } $F_{0}(\tau_{1},\tau_{2}) \not \in \{\eta_{k}:k<n(*)\}$,
moreover $F_{0}(\tau_{1},\tau_{2})\restriction 
\lg(\eta\conc \langle \rho \rangle)
 \not \in \{\eta_{k}:k<n(*)\}$,
\item[(7)] if $\eta\conc \langle \rho \rangle \in P^{5} \cap 
\{\eta_{k}:k<n(*)\}$, $\alpha<\beta<\kappa$,
$\tau_{1}=\eta\conc \langle \rho \rangle\conc\langle i_{1} \rangle =
\eta_{\alpha,k_{1}},$ $k_{1}\in [n(*),n)$,
$\tau_{2}=\eta\conc \langle \rho \rangle\conc\langle i_{2} \rangle =
\eta_{\beta,k_{2}},$ $k_{2}\in [n(*),n)$, $i_{1}\not =i_{2}$, {\em then}
$F_{l}(\tau_{1},\tau_{2})\restriction 
\lg(\eta\conc \langle \rho \rangle)
 \not \in \{\eta_{\alpha,k}:k<n, \alpha<\kappa\}$ for $l=0,1$.
\end{enumerate}

It is easy to satisfy $(4)-(6)$. To satisfy $(7)$ note that the
function $F(*,*)\restriction \lg(\eta\conc\rho)$ is $1-1$ by
2.2(4). Therefore we can choose a required sequence by induction of 
length $\kappa$.

Now we will show that if $\alpha_{0}, \ldots ,\alpha_{m(*)}<\kappa$,
then $\cal{B}_{T} \vDash \bigcap_{m\leq m(*)}a_{m} >0$.
It is enough to define $f:T\rightarrow \{0,1\}$ such that 
$(*)_{f,T}$ holds and 

\begin{equation*}\cal{B}_{0}\vDash t(f(\eta_{\alpha_{m},0}), \ldots ,
f(\eta_{\alpha_{m},n-1}))=1 \mbox{ \rm for } m\leq m(*).
\end{equation*}

Define $f(\nu)=1$ if and only if 
one of the following occurs:
\begin{enumerate}
\item[(a)] $\nu\in \{\eta_{\alpha_{m},k}:k<n,m\leq 
m(*),\}$ and $f_{\alpha_{m}}
(\nu)=1$.
\item[(b)] For some $\eta\conc \langle \rho \rangle \in P^{5} \cap \{
\eta_{k}:k<n(*)\}$ and $m_{1}\not =m_{2}$, and $k_{1},k_{2}
\in [n(*),n)$, we have $\eta\conc \langle \rho \rangle\conc\langle
i_{1}\rangle =\eta_{\alpha_{m_{1}},k_{1}}$,
$\eta\conc \langle \rho \rangle\conc\langle
i_{2}\rangle =\eta_{\alpha_{m_{2}},k_{2}}$
and  $\nu=F_{0}(\eta_{{\alpha_{m_{1}},k_{1}}},\eta_{\alpha_{m_{2}},k_{2}})$.
\end{enumerate}

\begin{lemma} \begin{enumerate} 
\item[(1)] $\cal{B}_{0}\vDash t(f(\eta_{\alpha_{m},0}), \ldots ,
f(\eta_{\alpha_{m},n-1}))=1$ for $m\leq m(*)$.
\item[(2)] $(*)_{f,T}$ holds.
\end{enumerate}
\end{lemma}

\Proof  (1) It suffices to prove that $f\restriction \{\eta_{\alpha_{m},0},
\ldots , \eta_{\alpha_{m},n-1}\} \subseteq f_{\alpha_{m}}$.
Assume first 
that $f_{\alpha_{m}}(\eta_{\alpha_{m},k})=1$. By the definition of 
$f$ we have $f(\eta_{\alpha_{m},k})=1$. Now assume that $f
(\eta_{\alpha_{m},k})=1$. Hence one of the cases (a) or (b) holds.
If case (a) holds we are done. So suppose that (b) holds. By (7)
it follows that $\eta_{\alpha_{m},k} \not \in \{\eta_{\alpha,k}:
\alpha<\kappa, k<n\}$, a contradiction.

\smallskip
\noindent (2) Assume that $(*)_{f,T}$
fails. Then there is ${\bold{e}}_{\nu_{1},\nu_{2}} \in \Gamma_{T}$
such that $\nu_{1}=\eta\conc\langle \rho \rangle \conc \langle i_{1}\rangle$.
$\nu_{2}=\eta\conc\langle \rho \rangle \conc \langle i_{2}\rangle$,
$i_{1}\not =i_{2}$ and $f(\nu_{1})=1=f(\nu_{2})$ and 
$f(F_{0}(\nu_{1},\nu_{2}))=f(F_{1}(\nu_{1},\nu_{2}))$.
Working toward a contradiction we consider three cases.

\smallskip
\noindent {\em Case 1.} $\nu_{1},\nu_{2} \in \{\eta_{\alpha_{m},k}:
m\leq m(*),
k<n\}$. Hence there is $m_{1},m_{2}\leq m(*)$, end 
$k_{1},k_{2}<n$ such that 
$\nu_{1}=\eta_{\alpha_{m_{1}},k_{1}}$, 
$\nu_{2}=\eta_{\alpha_{m_{2}},k_{2}}$.

If $m_{1}=m_{2}$, then as $\{\eta_{\alpha_{m_{1}},k}:k<n\}$ is 
$1$-closed, we have $\nu_{1},\nu_{2},F_{0}(\nu_{1},\nu_{2}),$
$F_{1}(\nu_{1},\nu_{2}) \in \{\eta_{\alpha_{m_{1}},k}:k<n\}$.
Since $(*)_{f_{\alpha_{m_{1}},T}}$ holds  we get a 
contradiction.

Hence $m_{1}\not =m_{2}$, and  $k_{1},k_{2}\in [n(*),n)$.
By the definition of $f$ we have 

\noindent $f(F_{0}(\nu_{1},\nu_{2}))=1$, 
so it suffices to show that $f(F_{1}(\nu_{1},\nu_{2}))=0$.
Assume to the contrary that $f(F_{1}(\nu_{1},\nu_{2}))=1$. Hence,
as $F_{1}(\nu_{1},\nu_{2}) \not \in \{\eta_{\alpha,k}:
\alpha<\kappa, k<n\}$ by (7), case (b) must hold. So there is 
${\eta^{*}}\conc\langle \rho^{*}\rangle \in P^{5}\cap \{\eta_{k}:k<n(*)\}$ and
$\tau_{1},\tau_{2}$ such that 
\begin{enumerate}
\item[(i)] $\tau_{1}={\eta^{*}}\conc\langle \rho^{*}\rangle\conc
\langle i_{3}\rangle=\eta_{\alpha_{m_{3}},k_{3}}$,
\item[(ii)] $\tau_{2}={\eta^{*}}\conc\langle \rho^{*}\rangle\conc
\langle i_{4}\rangle=\eta_{\alpha_{m_{4}},k_{4}}$,
\item[(iii)] $k_{3},k_{4}\in [n(*),n)$,
\item[(iv)] $m_{3}\not =m_{4}$,
\item[(v)] $F_{1}(\nu_{1},\nu_{2})=F_{0}(\tau_{1},\tau_{2})$.
\end{enumerate}

Note that ${\eta^{*}}\conc\langle \rho^{*}\rangle\not =
\eta\conc\langle \rho\rangle$. On the other hand (5) implies that 
$\eta=\eta^{*}$. Let $\zeta=\sup\{\xi:\rho^{*}\restriction \xi=
\rho\restriction \xi\}$. Since 
$\{\eta_{k}:k<n(*)\}$ is $2$-closed, it follows that 
${\eta}\conc\langle \rho^{*}\restriction \zeta \rangle =
\eta\conc\langle \rho\restriction \zeta \rangle \in \{\eta_{k}:k<n(*)\}$.
And moreover both ${\eta}\conc\langle \rho^{*}\restriction (\zeta
+1) \rangle$ and   
$\eta\conc\langle \rho\restriction (\zeta
+1) \rangle$ are in $\{\eta_{k}:k<n(*)\}$.


Note that $\lg(F_{1}(\tau_{1},\tau_{2})(\lg(\eta\conc\langle \rho\rangle)))
=\zeta$. 
Otherwise  
$\{F_{0}(\tau_{1},\tau_{2}), F_{1}(\tau_{1},\tau_{2},)\} \cap 
\{F_{0}(\nu_{1},\nu_{2}), F_{1}(\nu_{1},\nu_{2},)\}=\emptyset$.
Hence 
$F_{0}(\tau_{1},\tau_{2})\restriction \lg(\eta\conc\langle\rho\rangle)=
\eta\conc\langle \rho\restriction \zeta\rangle\in 
\{\eta_{k}:k<n(*)\}$,  contradicting $(7)$.

\medskip
\noindent {\em Case 2.}
$\{\nu_{1},\nu_{2}\} \cap  \{\eta_{\alpha_{m},k}:
m\leq m(*), k<n\}$ is a singleton.
By symmetry assume that $\nu_{1} \in \{\eta_{\alpha_{m},k}:
m\leq m(*), k<n\}$. Hence 
there are 
${\eta^{*}}\conc\langle \rho^{*}\rangle \in P^{5}\cap \{\eta_{k}:k<n(*)\}$,
$\tau_{1},\tau_{2}$ and  
$m\leq m(*)$, $i_{3} \not =i_{4}$, $k,k_{3},k_{4}<n$ such that
\begin{enumerate}
\item[(a)] $\nu_{1}=\eta_{\alpha_{m},k}$,
\item[(b)] $\nu_{2}=F_{0}(\tau_{1},\tau_{2})$,
\item[(c)]  $\tau_{1}={\eta^{*}}\conc\langle \rho^{*}\rangle\conc
\langle i_{3}\rangle$,
\item[(d)] $\tau_{2}={\eta^{*}}\conc\langle \rho^{*}\rangle\conc
\langle i_{4}\rangle$.
\end{enumerate}

It follows that ${\eta^{*}}\conc\langle \rho^{*}\restriction \zeta \rangle 
\vartriangleleft \nu_{2}=\eta\conc\langle \rho\rangle\conc\langle 
i_{2}\rangle$,
where $\zeta=h_{\sigma(\eta^{*})}(\tau_{1},\tau_{2})$.

As the last element of the sequence 
${\eta^{*}}\conc\langle \rho^{*}\restriction \zeta \rangle$ has length 
$\zeta$ 
we must have:
${\eta^{*}}\conc\langle \rho^{*}\restriction \zeta \rangle 
\vartriangleleft \eta$. 
Hence ${\eta^{*}}\conc\langle \rho^{*}\restriction \zeta \rangle\in 
\{\eta_{\alpha_{m},k}:k<n\}$. Since $\{\eta_{\alpha_{m},k}:k<n\}$ is 
$2$-closed and it follows by 2.3(2)(e) that $\nu_{2}\in 
\{\eta_{\alpha_{m},k}:k<n\}$, contradiction.

\medskip
\noindent {\em Case 3.}
$\{\nu_{1},\nu_{2}\} \cap  \{\eta_{\alpha_{m},k}:
m<m(*), k<n\}=\emptyset$.
By the definition there are $\tau_{l}$,
$m_{l}$, $k_{l}$  $i_{l}$ for $l=1,\ldots ,4$ 
such that $\nu_{1}=F_{0}(\tau_{1},\tau_{2})$, and 
$\nu_{2}=F_{0}(\tau_{3},\tau_{4})$, and 
\begin{enumerate}
\item[(a)] $\tau_{1}=\eta_{1}\conc\langle \rho_{1}\rangle\conc\langle
i_{1}\rangle =\eta_{\alpha_{m_{1},k_{1}}}$,
$\tau_{2}=\eta_{1}\conc\langle \rho_{1}\rangle\conc\langle
i_{2}\rangle =\eta_{\alpha_{m_{2},k_{2}}}$,
\item[(b)] $m_{1}\not =m_{2}$, $i_{1}\not =i_{2}$,
\item[(c)] $\tau_{3}=\eta_{2}\conc\langle \rho_{2}\rangle\conc\langle
i_{3}\rangle =\eta_{\alpha_{m_{3},k_{3}}}$,
$\tau_{4}=\eta_{2}\conc\langle \rho_{2}\rangle\conc\langle
i_{4}\rangle =\eta_{\alpha_{m_{4},k_{4}}}$,
\item[(d)] $m_{3}\not =m_{4}$, $i_{3}\not =i_{4}$, 
\item[(e)] $\eta_{1}\conc\langle \rho_{1}\rangle,
\eta_{2}\conc\langle \rho_{2}\rangle \in P^{5} \cap \{\eta_{k}:k\leq 
n(*)\}$.
\end{enumerate}

Let $\zeta_{1}=h_{\sigma(\eta_{1})}(\tau_{1},\tau_{2})$,
$\zeta_{2}=h_{\sigma(\eta_{2})}(\tau_{3},\tau_{4})$.
Note that $\eta_{1}\conc\langle\rho_{1}\restriction 
\zeta_{1}\rangle \vartriangleleft
\eta\conc\langle \rho\rangle$
and $\eta_{2}\conc\langle\rho_{2}\restriction 
\zeta_{2}\rangle \vartriangleleft
\eta\conc\langle \rho\rangle$.
Hence either $\eta_{1}\conc\langle\rho_{1}\restriction 
\zeta_{1}\rangle \vartriangleleft
\eta_{2}\conc\langle\rho_{2}\restriction 
\zeta_{2}\rangle$ or
$\eta_{2}\conc\langle\rho_{2}\restriction 
\zeta_{2}\rangle \vartriangleleft
\eta_{1}\conc\langle\rho_{1}\restriction 
\zeta_{1}\rangle$.
Assume the first case, the other is 
symmetric.
If $\eta_{1}\conc\langle\rho_{1}\restriction 
\zeta_{1}\rangle =
\eta_{2}\conc\langle\rho_{2}\restriction 
\zeta_{2}\rangle$, then $\zeta_{1}=\zeta_{2}=\sup\{\xi:
\rho_{1}\restriction \xi =\rho_{2}\restriction \xi$\}, as 
$\nu_{1}\not =\nu_{2}$. Hence both $\eta_{1}\conc\langle\rho_{1}\restriction 
\zeta_{1}\rangle$ and $\eta_{2}\conc\langle\rho_{2}\restriction 
\zeta_{2}\rangle$ are in $\{\eta_{k}:k<n(*)\}$, and by 2.3
$\nu_{1},\nu_{2}\in\{\eta_{k}:k<n(*)\}$, contradiction.
Therefore $\eta_{1}\conc\langle\rho_{1}\restriction 
\zeta_{1}\rangle\vartriangleleft\eta_{2}$.
Since $\eta_{2}\in\{\eta_{k}:k<n(*)\}$, also $\eta_{1}\conc\langle\rho_{1}
\restriction 
\zeta_{1}\rangle\in\{\eta_{k}:k<n(*)\}$. As above this implies that 
$\nu_{1}\in \{\eta_{k}:k<n(*)\}$, contradiction.

\section{The exact density}
In this section we prove that the topological density of the algebra 
$\cal{B}$ is $\mu$.

\begin{lemma}
$d(\cal{B}_{T})\geq \mu$.
\end{lemma}

\Proof
Assume that $d(\cal{B}_{T}) < \mu$ Hence there is a sequence
$\cal{D}^{\otimes}=\langle 
D^{\otimes}_{j}:j<d(\cal{B}_{T})\rangle$ of ultrafilters of
$\cal{B}_{T}$ such that for every $a\in \cal{B}_{T}\setminus \{0\}$ there is 
$j$ such that $a\in D^{\otimes}_{j}$.
Let 
\begin{align*}\cal{F}&=\{\sigma: \sigma<\mu \mbox{ \rm and there are }
\eta\mbox{ \rm and } \bar{\cal{D}} \mbox{ \rm such that: } \\
\begin{split}
&(a)\quad \eta\in T\setminus P^{7}, \sigma(\eta)\geq \sigma \\
&(b)\quad\bar{\cal{D}}=\langle D_{j}:j<\sigma\rangle \text{ is a sequence of 
ultrafilters of }\cal{B}_{T}, \\
&(c)\quad \text{ if }\eta \vartriangleleft \nu_{1}, 
\eta \vartriangleleft \nu_{2}, \nu_{1}\not =\nu_{2},
\text{ then } x_{\nu_{1}} 
\vartriangle x_{\nu_{2}} \in \bigcup_{j<\sigma} D_{j}\}
\end{split}
\end{align*}

\noindent Note that $\cal{F}\not =\emptyset$.  In particular 
$\langle \rangle$, 
$\cal{D}^{\otimes}$ witness that $d(\cal{B}_{T})\in \cal{F}$.
Let $\sigma^{*}=\min(\cal{F})$ and let $\eta^{*},\cal{D}^{*}$ witness
that $\sigma^{*}\in\cal{F}$. Note that $\sigma^{*}\geq \aleph_{0}$.
 Without loss of generality $\sigma^{*}=\sigma(\eta^{*})$,
(otherwise use ${\eta^{*}}\conc\langle\sigma^{*} \rangle$ instead of
$\sigma^{*}$.)

Now we choose by induction on $\zeta < (\sigma^{*})^{+}$, a sequence 
$\rho_{\zeta}$ such that:
\begin{enumerate}
\item[(1)] ${\eta^{*}}\conc \langle \rho_{\zeta} \rangle \in T$,
\item[(2)] $\lg(\rho_{\zeta})=1+\zeta$,
\item[(3)] $\xi < \zeta \implies \rho_{\xi} =\rho_{\zeta}\restriction
(1+\xi)$,
\item[(4)] ${\eta^{*}}\conc \langle \rho_{\zeta+1} \rangle\in P^{6}$,
$\rho_{\zeta +1}(\zeta)=(\nu_{0,\zeta},\nu_{1,\zeta})$,
\item[(5)] if $\zeta=\sigma^{*}\xi+j$, $j<\sigma^{*}$,
then $(x_{\nu_{0,\zeta}} \vartriangle x_{\nu_{1,\zeta}}) \not \in
\bigcup_{i<j}D^{*}_{i}$
\end{enumerate}

For $\zeta=0$, let $\rho_{0}=\langle\sigma^{*}\rangle$.
For $\zeta$ limit put $\rho_{\zeta}=\bigcup_{\xi<\zeta}\rho_{\xi}$.
For $\zeta=\xi+1$, $\zeta=\sigma^{*}\xi+j$, if 
we cannot find suitable $(\nu_{0},\nu_{1})$, then
${\eta^{*}}\conc \langle \rho_{\xi}\rangle$, $\langle D^{*}_{i}:i<j\rangle$
witness that $j\in\cal{F}$, contradicting the minimality of $\sigma^{*}$.

Let $\rho^{*}=\rho_{(\sigma(\eta^{*}) )^{+}}$.
Note that for every $i_{1}<i_{2}<\sigma^{*}$, we have 
$x_{{\eta^{*}}\conc \langle \rho^{*}\rangle\conc\langle i_{1}\rangle}
\vartriangle
x_{{\eta^{*}}\conc \langle \rho^{*}\rangle\conc\langle i_{2}\rangle}
\in \bigcup_{j<\sigma^{*}}D^{*}_{j}$. 
Hence there is $j(*)<\sigma^{*}$ such that 
$X=\{i<(\sigma^{*})^{+}:x_{{\eta^{*}}\conc \langle 
{\rho^{*}}\rangle\conc\langle i\rangle}\in D^{*}_{j(*)}\}$ has 
cardinality 
$(\sigma^{*})^{+}$.
Now let $j_{1}\in (j(*),\sigma^{*})$ and $i_{0},i_{1}\in X$ be such that 
$h_{\sigma^{*}}(\{i_{0},i_{1}\})=j_{1} \mod \;\sigma^{*}$, 
($j_{1} $ exists be the 
definition of $h_{\sigma^{*}}$).

Now let $\tau_{0}={\eta^{*}}\conc 
\langle \rho^{*}\rangle\conc\langle i_{0}\rangle$,
$\tau_{1}={\eta^{*}}\conc \langle \rho^{*}\rangle\conc\langle i_{1}\rangle$.
Hence $F_{0}(\tau_{0},\tau_{1})=\nu_{0,\epsilon}$ and 

\noindent  $F_{1}(\tau_{0},\tau_{1})=\nu_{1,\epsilon}$, where $\epsilon=
h_{\sigma^{*}}(\{i_{0},i_{1}\})$.
So ${\bold{e}}_{\tau_{0},\tau_{1}} \in \Gamma_{T}$.
Note that $x_{\tau_{0}}\cap x_{\tau_{1}}\in D^{*}_{j(*)}$, and 
$x_{\nu_{0,\epsilon}}\vartriangle x_{\nu_{1,\epsilon}} \not \in 
D^{*}_{j_{1}}$, as $j(*)<j_{1}$. Hence
$-(x_{\nu_{0,\epsilon}}\vartriangle x_{\nu_{1,\epsilon}}) \in 
D^{*}_{j(*)}$.
On the other hand 
${\bold{e}}_{\tau_{0},\tau_{1}} \in \Gamma_{T}$ implies that 
$x_{\tau_{0}}\cap x_{\tau_{1}}\cap
-(x_{\nu_{0,\epsilon}}\vartriangle x_{\nu_{1,\epsilon}})=0$,
contradiction.

\begin{lemma}
$d(\cal{B}_{T})\leq \mu$.
\end{lemma}

\Proof
The idea of the proof is to define a set $\cal{F} \subseteq
\{1,-1\}^{T}$ such that:

\begin{enumerate}
\item[$(1)$] $|\cal{F}|=\mu$
\item[$(2)$] For every $f\in \cal{F}$, the set 
$D_{f}=\{x_{\eta}^{f(\eta)}:\eta \in T\}$ generates an ultrafilter in 
$\cal{B}_{T}$, where $x^{1}=x$ and  $x^{-1}$ is the complement of $x$. 
\item[$(3)$] For every $a \in \cal{B}_{T}$, non-zero,
there is $f\in \cal{F}$ such that $a$ is an the ultrafilter generated by 
$D_{f}$.
\end{enumerate}

First, divide $T$ into three disjoint sets $T=T_{0}\cup T_{1}
\cup T_{2}$ as follows. 
$T_{0}=\bigcup \{\{\tau_{0},\tau_{1}\}:\{\tau_{0},\tau_{1}\} \in 
\dom(F^{i}), i=0,1\}$.
Let $T_{1}$ be the image of $T_{0}$ under $F^{0}$ and $F^{1}$.  
It follows from the construction of $T$ that $T_{0}$ is disjoint from 
$T_{1}$.
Finally let $T_{2}=T\setminus (T_{0} \cup T_{1})$.

Let $F_{2} \subseteq \{1,-1\}^{T_{2}}$ be a set of cardinality $\mu$ such that 

\begin{enumerate}
\item[$(a)$] for every finite $u\subseteq T_{2}$ and a function
$h:u \rightarrow \{1,-1\}$ there is $f\in F_{2}$ such that 
$h \subseteq f$.
\end{enumerate}

Similarly, let $F_{1}\subseteq \{1,-1\}^{T_{1}}$ be a set of cardinality $ \mu$
such that 

\begin{enumerate}
\item[$(b)$] for every finite $u\subseteq T_{1}$ and a function
$h:u \rightarrow \{1,-1\}$ there is $f\in F_{1}$ such that 
$h \subseteq f$.
\end{enumerate} 

Now, for every $f\in F_{1}$ we define a set $F_{0}^{f} \subseteq
\{1,-1\}^{T_{0}}$ of cardinality $\mu$ such that 

\begin{enumerate}
\item[$(c)$] for every $g\in F_{0}^{f}$, for every 
$\nu_{0},\nu_{1} \in T_{1}$, $\tau_{0},\tau_{1}\in T_{0}$,
if $F^{i}(\tau_{0},\tau_{1})=
\nu_{i}$, $i=0,1$ and $f(\nu_{0})=f(\nu_{1})$, then
$(g(\tau_{0}),g(\tau_{1})) \neq (1,1)$
\item[$(d)$] $F_{0}^{f}$ is dense with respect to $(c)$, i.e., 
for every finite $u \subseteq T_{0}$ and a function 
$h:u \rightarrow \{1,-1\}$ such that the condition $(c)$ is satisfied
with $h$ in place of $g$, then there is $g\in F_{0}^{f}$ such that
$h\subseteq g$, i.e. $(*)_{f \cup g, T}$ holds.
\end{enumerate}

Finally define $\cal{F}=\{f_{2} \cup f_{1} \cup f_{0}:
f_{2} \in F_{2},  f_{1}\in F_{1}, f_{0} \in F_{0}^{f_{0}}\}$.

We prove that $\cal{F}$ is as required. 
It is obvious that $(1)$ holds, and (2) follows from the definition of 
$\cal{F}$, $D_{f}$ and Proposition 3.2.

To prove $(3)$ let $\cal{B}_{T} \models a > 0$. Without loss of generality
$a= \bigcup_{\eta \in u} x_{\eta}^{h(\eta)}$ for some finite 
$u \subseteq T$ and $h:u \rightarrow \{1,-1\}$. Let $u_{i}=u \cap
T_{i}$ for $i \leq 2$. Let $f_{2} \in F_{2}$ and $f_{1} \in F_{1}$
be such that $h\restriction u_{i} \subseteq f_{i}$, $i=1,2$.
Let $f_{0} \in F^{f_{1}}_{0}$ be such that $h\restriction u_{0}
\subseteq f_{0}$. We have to show that $f_{0}$ exists. Note that   
$(*)_{h,T}$ holds
since $a$ is non-zero, hence $f_{0}$ exists by $(d)$.
It follows that $a \in D_{f}$, where $f=f_{2} \cup f_{1} \cup f_{0}$.

This finishes the proof
of the theorem.

\end{document}